\documentclass[12pt]{amsart}

\usepackage{graphicx}
\usepackage{amssymb}
\usepackage{array}
\usepackage{amsthm}
\usepackage{mathtools}
\usepackage[hidelinks]{hyperref}
\usepackage{tgtermes}
\usepackage{mathrsfs}
\usepackage{stmaryrd}
\usepackage[all]{xy}
\usepackage[mathcal]{eucal}
\usepackage{verbatim} 
\usepackage{fullpage}  
\usepackage{hyperref}
\usepackage{setspace}
\usepackage{tikz-cd}
\usepackage{amsmath}
\usetikzlibrary{positioning}
\DeclareMathOperator{\GL}{GL}

\DeclareMathOperator{\SL}{SL}

\newtheorem{theorem}{Theorem}[section]

\theoremstyle{definition}

\newtheorem{example}[theorem]{Example}
\newtheorem{remark}[theorem]{Remark}

\title{On The Klein--Hilbert Resolvent Problem \\ \small (\"{U}ber Das Klein--Hilbertsche Resolventenproblem)}
\author{G.N. Chebotarev \\Translator: Sidhanth Raman}
\date{August 2023}

\begin{document}

\begin{abstract}
This is an English translation of G.N. Chebotarev's paper ``On The Klein--Hilbert Resolvent Problem" which was originally written in Russian and published in Izvestiya Kazan. Fiz. Mat. Obshch., 6 (1932–1933), 5–22. In this article, Chebotarev addresses ``Klein's resolvent problem"; in modern parlance, Chebotarev gives upper bounds on the essential dimension of finite simple groups via their representation theory.
\end{abstract}

\maketitle

Some years ago I published two papers (\cite{tschebotarow1} and \cite{tschebotarow2}), in which I set forth my results relating to an algebraic problem posed by F. Klein under the name ``problem of forms" \cite{klein1922gesammele} and subsequently greatly generalized by D. Hilbert (\cite{hilbert1900gott} and \cite{hilbert1927gleichung}). However, since my works were written \textit{in statu nascendi}, so that not everything was clear to me during the editing, while now I am able to overcome some of the difficulties encountered there, I believe it is appropriate to give a new systematic exposition of my results. In doing so, I mostly restrict myself to the simplest case, when the Galois group is simple, although some of my theorems easily extend to more general cases as well. However, the case considered here is almost the only essential one for this theory. This article does not presume any prior familiarity with my earlier works. 

I will first present the two algebraic problems this article focuses on, and then I will explain the difference between them. The first of which, basically posed by Klein, can be considered a special case of the second problem posed by Hilbert.\\

\noindent
I. \textbf{Klein's problem}. We are given an algebraic equation
\begin{align}\label{eq1}
    f(x) = x^n + a_1x^{n-1} + \dots + a_{n-1}x + a_n = 0,
\end{align}
whose coefficients $a_1,\dots,a_n$ are assumed to be arbitrary variables. Our goal is to find a Tschirnhaus transformation of this equation such that the coefficients of the transformed equation would depend on the least possible number of variable parameters. 

The coefficients of the sought transformation may contain numerical irrationalities and some function $\Phi(x_1, x_2, \dots, x_n)$ of the roots $x_1, x_2, \dots, x_n$ of Equation \ref{eq1} belonging to a given group of substitutions $\mathcal{G}$.

This problem can be formulated in terms of field theory as follows:
\begin{enumerate}
    \item[(Ia)] We are given a field $k$ containing $n+1$ variables $a_1, a_2, \dots, a_n; \Phi$, which are related by an algebraic relation
\begin{equation*}
F(a_1, a_2, \dots, a_n; \Phi)=0.
\end{equation*}
In addition, an algebraic extension $K/k$ of this field is given, whose generating elements $x_1, x_2, \dots, x_n$ are related to the generating elements of the field $k$ as follows:
\begin{align*}
x_1 + x_2 + \dots + x_n &= -a_1, \\  
x_1x_2 + \dots + x_{n-1}x_n &= a_2, \\
&\vdots \\
x_1x_2\dots x_{n-1} + \dots + x_2x_3\dots x_n &= (-1)^{n-1}a_{n-1}, \\
x_1x_2\dots x_n &= (-1)^na_n, \\
\Phi &= \Phi(x_1, x_2, \dots, x_n),
\end{align*}
where $\Phi(x_1, x_2, \dots, x_n)$ is a rational function belonging to the given group of substitutions $\mathcal{G}$. It is assumed here that Equation \ref{eq1} is identically satisfied after the above substitutions. It is required to find such a subfield $\mathfrak{F}$ of the field $K$ that:
\begin{enumerate}
     \item[(1)] the compositum of the fields $\mathfrak{F}$ and $k$ is the field $K$, and
     \item[(2)] the degree of transcendence $s$ of the field $\mathfrak{F}$ is as small as possible.
\end{enumerate}
The number $s$ will be called the true degree of transcendence of the extension $K/k$.\footnote{Translator's note: In modern language, this number $s$ in the ``Klein problem" is called the \textit{essential dimension} of the extension $K/k$, denoted $\textrm{ed}(K/k)$.}
\end{enumerate}
The latter formulation of Klein's resolvent problem admits the following natural generalization:
\begin{enumerate}
\item[(Ib)] Given is an arbitrary field $k$ of transcendence degree $n$ and its algebraic extension $K/k$. Find a subfield $\mathfrak{F}$ of the field $K$ with the least possible degree of transcendence $s$ such that the fields $\mathfrak{F}$ and $k$ generate the entire field $K$.
\end{enumerate}

\noindent
II. \textbf{Hilbert's resolvent problem}. Suppose we are given an algebraic extension $K/k$ of the field $k$ with transcendence degree $n$. The goal is to find the smallest integer $s$ possessing the following properties\footnote{Translator's note: In modern language, this number $s$ in the ``Hilbert problem" is called the \textit{resolvent degree} of the extension $K/k$, denoted $\mathrm{RD}(K/k)$.}:
\begin{enumerate}
    \item there exists a sequence of algebraically-transcendental extensions $$K_1/k, K_2/K_1, \dots, K_m/K_{m-1},$$ each with true degree of transcendence equal to $s$; 
    \item the field $K_m$ contains the field $K$.
\end{enumerate}

\noindent
In this work, I completely solve only Problem Ia. It seems very likely that Problem II is not a trivial generalization of problem I, i.e. that in some cases the solution of Problem II gives a smaller value for $s$ than the solution of Problem I.\footnote{Translator's note: This is indeed the case --- there are instances in which the resolvent degree of an extension $K/k$ is strictly bounded above by its essential dimension. Note that lower bounds on essential dimension have been computed, but to date, there are no nontrivial lower bounds on resolvent degree.}\\

\tableofcontents

\section{The theory of representations of continuous transformation groups}

Let us consider continuous transformation groups\footnote{Translator's note: In modern parlance, a ctg is defined to be not just a Lie group, but a Lie group along with a (linear) representation.} (for short: ctg) of the form
\begin{equation}\label{eq2}
x_i' = f_i(x_1, x_2, \dots, x_n; a_1, a_2, \dots, a_r) \quad (i = 1, 2, \dots, n), \tag{1.1}  
\end{equation}
where two ctg's are called \textit{locally isomorphic} (according to Lie: \textit{gleich zusammengesetz}) if the infinitesimal operators $X_1, X_2, \dots, X_r$ of both groups can be chosen so that the structural constants $c^k_{ij}$ in the relations 
\begin{equation*}
(X_i, X_j) = \sum_{k=1}^r c^k_{ij} X_k \quad (i,j = 1, 2, \dots, r)
\end{equation*}
are the same for both groups.\footnote{Translator's note: This notion of local isomorphism of the ctg's in question is an isomorphism of their corresponding Lie algebras.}

Two ctg's $\Gamma$ and $\Gamma'$ are called \textit{similar} if there exists a transformation $S$ that maps each transformation $A$ of the group $\Gamma$ to some transformation $A'$ of the group $\Gamma'$: $S^{-1}AS = A'$ or $S^{-1}\Gamma S = \Gamma'$. More precisely, if
\begin{align*}
\Gamma: x' &= f(x,a), \\
\Gamma': y' &= f'(y,a),
\end{align*}
then the two ctg's are similar if there exists an invertible transformation
\begin{equation*}
S: \quad x=\varphi(y), \quad 
y=\varphi^{-1}(x),
\end{equation*}
such that
\begin{align*}
\varphi^{-1}(f'(\varphi(y),a)) &= f(y,a) \\
f'(\varphi(y),a) &= \varphi(f(y,a)).
\end{align*}

For simplicity, we assume that the parameter $a$ is the same in both groups. In general, to transform the group $\Gamma$ into $\Gamma'$, one must also transform the parameter $a$ along with the variable $x$. But groups that are transformed into each other by a transformation of the parameter will be considered as simply coinciding. 

Thus, in order to go from group $\Gamma$ to group $\Gamma'$, one must replace the variables $x$ in the equations defining the group with the variables $y$ using the transformation $S$.

We also consider the case when the transformation $S$ is irreversible, i.e. when the system of functions $\varphi_i$ contains fewer than $n$ functionally independent functions. In this case, we say that the representation $\Gamma'$ is contained in the representation $\Gamma$, or symbolically: $\Gamma' \preceq \Gamma$. Clearly, this concept is transitive, i.e.
\begin{theorem}
If $\Gamma_1 \preceq \Gamma$ and $\Gamma_2 \preceq \Gamma_1$, then $\Gamma_2 \preceq \Gamma$.
\end{theorem}

Let $\Gamma$ be the parametric group of some ctg $\mathcal{G}$, and $\Gamma_1$ an arbitrary representation of the group \ref{eq2}. Then
\begin{theorem}\label{thm2} 
The parametric group of a ctg contains all existing representations of this group.
\end{theorem}

\begin{proof}
Let \ref{eq2} be the equations of the ctg $\Gamma_1$, and 
\begin{equation} 
a_i = \varphi_i(a,b) \quad (i=1,2,\dots,r) \tag{1.2}
\end{equation}
the equations of its parametric group $\Gamma$. Then
\begin{equation}\label{ctg2}
f_i(f(x,a),b) = f_i(x,\varphi(a,b)) \quad (i=1,2,\dots,n). \tag{1.3}
\end{equation}
Assigning arbitrary constant values $x^{(0)}$ to the variables $x$, equations \ref{ctg2} show that if we subject the parameters $a_i$ to the transformations $a_i = \varphi_i(a,b)$ of the group $\Gamma$, then the variables $x_i' = f_i(x^{(0)}, a)$ will undergo the transformations $x_i'' = f_i(x',b)$ of the group $\Gamma$, and so on. 
\end{proof}

A key characteristic of the parametric group is its property of being \textit{simply transitive}. This means that in it there is one and only one transformation that maps one of two arbitrarily given points to the other. Analytically, this property is expressed by the fact that the order of the group and the dimensionality of its representation space coincide and that its infinitesimal operators are not related by any nontrivial linear dependence of the form
\begin{equation*}
\varphi_1(x_1, x_2, \dots, x_n) X_1(f) + \varphi_2(x_1, x_2, \dots, x_n) X_2(f) + \dots + \varphi_n(x_1, x_2, \dots, x_n) X_n(f) = 0.
\end{equation*}

Theorem \ref{thm2} can be proven relying only on the simple transitivity of the parametric group. Indeed, let $X_1, X_2, \dots, X_r$ and $Y_1, Y_2,\dots,Y_r$ be systems of infinitesimal operators of two locally isomorphic groups, and let the first one be simply transitive. Then the system 
\begin{equation}\label{eq14}
    X_i(f) + Y_i(f) = 0 \quad (i = 1,2,\dots,r)
    \tag{1.4}
\end{equation}
has exactly $(n + r) - r = n$ independent integrals $\Phi_1,\Phi_2,\dots,\Phi_n$, where $n$ is the dimension of the representation space $(y_1, y_2, \dots, y_n)$ of the first group. By virtue of the transitivity of the first group, system \ref{eq14} has no solutions depending only on $X_i$, and therefore the Jacobian 
$$\frac{d(\Phi_1,\Phi_2,\dots,\Phi_n)}{d(y_1,y_2,\dots,y_n)}$$ 
does not vanish identically. Consequently, the system of finite equations 
$$\Phi_1 = c_1 , \Phi_2 = c_2 ,\dots , \Phi_n = c_n$$
is solvable with respect to $y_1, ..., y_n$ and it is easy to see that the equations obtained
\begin{equation}\label{eq15}
\begin{matrix}
    y_1 - \varphi_1(x_1,x_2\dots,x_n) = 0\\
y_2 - \varphi_2(x_1,x_2\dots,x_n) = 0 \\
\vdots\\
y_n - \varphi_n(x_1,x_2\dots,x_n) = 0
\end{matrix}
\tag{1.5}
\end{equation}
define the desired transformation \cite{lie1888theorie}. 

It follows directly from this that two simply transitive locally isomorphic ctg's are similar. 

Let $\Gamma$ and $\Gamma_1$ be two representations of the same ctg, with $\Gamma \preceq \Gamma_1$. The latter always holds if the group $\Gamma$ is simply transitive. If the substitution \ref{eq15} transforms $\Gamma$ into $\Gamma_1$, then the equations 
\begin{equation}\label{eq16}
\varphi_i(x_1,x_2,\dots,x_r) = c_i \quad (i=1,2,\dots,r) \tag{1.6}
\end{equation}
define one of the systems of imprimitivity of the group $\Gamma$ due to the fact that any transformation of the group $\Gamma$ transforms the variables $y_i$ into functions of the $y_i$ alone. Conversely, if the equations \ref{eq16} represent a system of imprimitivity of the group $\Gamma$, then each transformation of the group $\Gamma$ transforms the equations \ref{eq16} into equations of the form $\psi_i(x_1, x_2, \dots, x_n) = c_i'$, where the constants $c_i'$ are functions of the $c_i$ and parameters $a$ of the group $\Gamma$: 
$$c_i' = \psi_i(c',a) \quad (i=1,2,\dots,n).$$
In other words, we have the relations
$$\varphi_i(f(x,a)) = \psi_i(\varphi(x), a) \quad (i=1,2,\dots,n).$$ 
It follows from this that the equations 
$$y_i = \psi_i(y_i, a) \quad (i=1,2,\dots,n)$$
define an ctg isomorphism of $\Gamma$ to $\Gamma_1$, which is obtained from $\Gamma$ by means of the transformation 
$$y_i = \varphi_i(x_1, x_2, \dots, x_n).$$
This group is transitive if and only if the functions $\varphi_i(x)$ are functionally independent. Deleting from the system $\varphi_i(x)$ those functions that depend on the others, we obtain a transitive representation (a \textit{truncated group}). Conversely, when given a transitive representation, one can obtain from it some intransitive representation if one supplements the system of variables $y_i = \varphi_i(x_1, x_2, \dots, x_n)$ with some new system of variables invariant under the transformations of the group $\Gamma$. If one then performs the most general transformations over the system of variables obtained in this way, one gets the most general representation belonging to our system of imprimitivity. What has been said allows us to consider only transitive representations. 

If the representation $\Gamma$ contains another representation $\Gamma_1$ of the same group and the latter is not similar to the representation $\Gamma$, then the variables of the representation $\Gamma_1$ expressed as functions of the variables of the representation $\Gamma$ determine the left-hand sides of the equations of some system of imprimitivity of the group $\Gamma$. On the other hand, Lie (\cite[p. 522, Theorem 92]{lie1888theorie}) found the following most general method of constructing systems of imprimitivity of a transitive group. One should construct the stabilizer group (\textit{stabilit\"{a}tsgruppe}) $H_P$ of an arbitrary point $P$ of the representation space of $\Gamma$, i.e. the set of all transformations of the group $\Gamma$ that leave the point $P$ fixed. Then take an arbitrary group $\Sigma$ lying between $\Gamma$ and $H_P$ (the existence of such a group is necessary and sufficient for the group $\Gamma$ to be imprimitive) and determine the manifold obtained by subjecting the point $P$ to all transformations of the group $\Sigma$. The number of degrees of freedom of this system of imprimitivity (i.e. the degree of functional independence among the left-hand sides of the equations of the system of imprimitivity; in other words, the dimension of the representation space of the truncated group $\Gamma$) is equal to the index of the group $\Sigma$ relative to $\Gamma$, i.e. the difference between the orders of these groups. The representation $\Gamma_1$ will be proper (\textit{faithful}\footnote{Translator's note: The original text contained the German word \textit{tre\'{u}}. Tre\'{u} means faithful in the modern theory of groups and representations.}), i.e. isomorphic to $\Gamma$ if and only if $\Sigma$ does not contain any normal divisors\footnote{Translator's note: \textit{Divisor} is the term for subgroup, and so normal divisors are normal subgroups.} of the group $\Gamma$ other than the identity group.

If the group $\Gamma$ is simply transitive, then the stabilizer group $H$ coincides with the identity group, so that each subgroup of the parametric group corresponds to a contained representation of $\Gamma$, and vice versa. This representation is primitive if and only if the corresponding subgroup is maximal, i.e. is not contained in any subgroup of the group $\Gamma$ other than $\Gamma$ itself. If $\Gamma_1$ and $\Gamma_2$ are two representations of the same group, then $\Gamma_1 \preceq \Gamma_2$ if and only if for the corresponding subgroups $\Sigma_1, \Sigma_2$ of the parametric group we have 
$$\Sigma_1 \leq S^{-1}\Sigma_2 S,$$
where $S$ is a suitably chosen transformation of the parametric group, and the symbol $\leq$ means that the left-hand side is a divisor of the right-hand side.

We thus arrive at the following theorem:

\begin{theorem}
A ctg $\Gamma$ has a representation in an $s$-dimensional space (i.e. $\Gamma$ is an $s$-group) if and only if it has a subgroup of index $s$ not containing any normal divisor of the group $\Gamma$ other than the identity group.
\end{theorem}

It seems likely that any representation can be obtained by purely algebraic operations from the finite equations of the parametric group. We will first show that the infinitesimal operators of such a representation can be obtained algebraically, and then prove that, under certain restrictions imposed on the structure of the group, the same is possible for its finite equations. 

To prove the first part of this assertion, it is sufficient to show that all subgroups of a given transitive ctg can be defined algebraically. Lie (\cite[p. 208-210, Section 5]{lie1888theorie}) gave the following way to make this definition. Let $X_1, X_2, \dots, X_r$ be the infinitesimal operators of the ctg $\Gamma$. We seek such linear combinations of them
\begin{equation}\label{eq17}
Y_i = g_{i1}X_1 + g_{i2}X_2 + \dots + g_{ir}X_r \quad (i=1,2,\dots,m), \tag{1.7}
\end{equation}
that satisfy the relations
\begin{equation}\label{eq18}
(Y_i,Y_j) = \sum_{v=1}^m \gamma_{ij}^v Y_v \quad (i,j=1,2,\dots,m). \tag{1.8}
\end{equation}
Substituting in the latter $Y_i$ from \ref{eq17} and using the formulas
\begin{equation*}
(X_k, X_l) = \sum_{v=1}^r c^v_{kl} X_v \quad (k,l=1,2,\dots,r),
\end{equation*}
we get
$$(Y_i,Y_j) = \sum_{\lambda,\mu = 1}^r g_{i\lambda} g_{j \mu} (X_\lambda, X_\mu) = \sum_{\lambda,\mu,\nu} g_{i\lambda}g_{j\mu} c_{\lambda \mu}^\nu X_\nu.$$
On the other side of the equation, 
$$\sum_{\nu} \gamma_{ij}^\nu Y_\nu = \sum_{\nu, \mu} \gamma_{ij}^\nu g_{\nu\mu} X_{\mu}.$$
Eliminating the constants $\gamma_{ij}^\mu$ from this system, we obtain a system of algebraic equations for determining $g_{\mu\nu}$. Substituting the solution of this system into \ref{eq17} gives the sought infinitesimal operators $Y_1, Y_2,\dots, Y_m$

Before proceeding to the definition of the finite equations of the subgroup, we will prove one very important theorem of \'{E}. Cartan \cite{cartan1902structure}, which allows us to consider the set of all representations of a ctg as a set equipped with an ``Archimedean property" in a certain sense.

\begin{theorem}\label{thm4} 
If $\Gamma$ and $\Gamma_1$ are any two representations of a ctg, then there exists some degree $m$ such that $\Gamma^m \preceq \Gamma_1$.
\end{theorem}

By the degree $\Gamma^m$ representation of
\begin{equation*}
\Gamma: \quad x_i' = f_i(x_1, x_2, \dots, x_n ; a_1 , a_2, \dots , a_r) \quad (i = 1,2,\dots, m)
\end{equation*}
we mean a representation in an $mn$-dimensional space 
\begin{align*}
x_1, x_2, &\dots, x_n, \\
x_1', x_2', &\dots, x_n',\\
&\vdots \\
x_1^{(m-1)}, x_2^{(m-1)}, &\dots, x_n^{(m-1)}, 
\end{align*}
where any transformation induces a transformation of the representation $\Gamma$ with the same parameter values $a_1, a_2, \dots, a_r$ in each space $x_1^{(i)}, x_2^{(i)}, \dots, x_n^{(i)}$ (an ``extended group", according to Lie).

\begin{proof}
Let 
\begin{equation*}
X_i (f) = \sum_{\nu} \xi_{i\nu}(x) \frac{df}{dx_\nu} \quad (i=1,2,\dots,r)
\end{equation*}
be the infinitesimal operators of the representation $\Gamma$. Choose $m$ large enough so that the rank of the matrix
\begin{equation*}
\begin{matrix}
\xi_{11}(x), \xi_{12}(x), \dots, \xi_{1n}(x); \xi_{11}(x'), \dots ; \dots, \xi_{1n}(x^{(m-1)}), \\
\xi_{21}(x), \xi_{22}(x), \dots, \xi_{2n}(x); \xi_{21}(x'), \dots ; \dots, \xi_{2n}(x^{(m-1)}), \\
\vdots\\
\xi_{r1}(x), \xi_{r2}(x), \dots, \xi_{rn}(x); \xi_{r1}(x'), \dots ; \dots, \xi_{rn}(x^{(m-1)}) \\
\end{matrix}  
\end{equation*}
is equal to $r$. This is always possible since, due to the independence of the operators $X_1, X_2, \dots, X_r$, there are no linear relations with constant coefficients between the rows of this matrix. 

Then find the solutions of the system
\begin{equation*}
\begin{matrix}
Y_1(f) + X_1(f) + X_1'(f) + \cdots + X_1^{(m-1)}(f) &= 0, \\
Y_2(f) + X_2(f) + X_2'(f) + \cdots + X_2^{(m-1)}(f) &= 0, \\
\vdots \\
Y_r(f) + X_r(f) + X_r'(f) + \cdots + X_r^{(m-1)}(f) &= 0,
\end{matrix}
\end{equation*}
where 
$$X_j^{(i)}(f) = \sum_{\nu=1}^n \xi_{j\nu}(x^{(i)}) \frac{df}{dx^{(i)}},$$ 
and $Y_i(f)$ is the operator corresponding to the operator $X_i(f)$ in the representation $\Gamma_1$ given in the space $(y_1, y_2, \dots, y_s)$. Clearly, this system is complete and therefore has $s + mn - r$ independent integrals. Among these integrals there are exactly $mn - r$ that depend only on the variables $x_j^{(i)}$ ($i=0,1,\dots,m-1; j=1,2,\dots,n$). The remaining integrals can therefore be reduced to the form
\begin{equation*}
\begin{matrix}
y_1 - \theta_1(x_1, x_2, \dots, x_n^{(m-1)}), \\
y_2 - \theta_2(x_1, x_2, \dots, x_n^{(m-1)}), \\
\vdots \\  
y_s - \theta_s(x_1, x_2, \dots, x_n^{(m-1)}).
\end{matrix}
\end{equation*}
Then the formulas
\begin{equation*}
\begin{matrix}
y_1 = \theta_1(x_1, x_2, \dots, x_n^{(m-1)}), \\
y_2 = \theta_2(x_1, x_2, \dots, x_n^{(m-1)}), \\
\vdots \\
y_s = \theta_s(x_1, x_2, \dots, x_n^{(m-1)}) 
\end{matrix}
\end{equation*}
give the sought transformation that transforms the representation $\Gamma^m$ into the representation $\Gamma_1$.
\end{proof}

\begin{remark}
We do not assume anywhere that the representation $\Gamma_1$ is proper. Thus, one can state that a sufficiently high degree of any representation of a given ctg $\mathcal{G}$ also contains representations of the quotient groups of the group $\mathcal{G}$.\footnote{Translator's note: There are assumptions missing from this assertion; for example, no power of the representation $\det: \GL_n \to \GL_1$ contains a nontrivial representation of $\SL_n$.}
\end{remark}

We now turn to proving our assertion about the possibility of algebraically finding the finite representation equations. In doing so, we will assume that the group $\mathcal{G}$ is \textit{semisimple}. This means that it does not contain any solvable normal divisors other than the identity. Then its adjoint group $E$ is a proper representation. By virtue of Theorem \ref{thm4}, a sufficiently high degree of the group $E$ contains the parametric group, and therefore we can consider such a degree of the group $E$ instead of the parametric group. In order to algebraically obtain any maximal subgroup of the group $\mathcal{G}$, we assume, following W. Killing \cite{killing1890bestimmung}, that the sought subgroup contains a transformation whose nonzero characteristic roots are all simple. Cartan \cite{cartan1894structure} considers it highly plausible that this assumption does not carry any real restrictions. It is known that the infinitesimal operators of the group $E$ can be normalized in such a way that some of them, $Y_1, Y_2, \dots, Y_\ell$ ($\ell$ is called the \textit{rank} of the group $\mathcal{G}$), commute with each other, while the remaining $X_1, X_2, \dots, X_{r-\ell}$ are related to the $Y_i$ by the relations  
$$(X_i, e_1Y_1 + e_2Y_2 + \cdots + e_\ell V_\ell) = (e_1\omega_i^{(1)} + e_2\omega_i^{(2)} + \cdots + e_\ell\omega_i^{(\ell)})X_i,$$
where $\omega_i^{(j)}$ are the nonzero characteristic roots of the operator $Y_j$. 

Killing (\cite[p. 242]{killing1890bestimmung}) showed that any subgroup $\mathcal{H}$ of the group $\mathcal{G}$ containing, for example $Y_1$, is generated by the operators $X_i, Y_j$ contained in it. If, in addition, the group $\mathcal{H}$ is maximal, then it contains all the operators $Y_1, Y_2, \dots, Y_\ell$. One can take for $Y_1$ such a linear combination of $Y_1, Y_2, \dots, Y_\ell$ that the corresponding characteristic roots are simple and their ratios are rational (for this, it is sufficient to take as $Y_1$ one of the operators $(X_{\alpha}, X_{\alpha'})$, where $X_{\alpha}$ and $X_{\alpha'}$ are operators corresponding to the roots $\omega_{\alpha}$ and $\omega_{\alpha'}$ (\cite[p. 42]{cartan1894structure}). Moreover, the whole abelian group generated by the operators $Y_i$ can be constructed from operators of this kind. On the other hand, it is obvious that the whole group $\mathcal{G}$ is generated by means of operators whose nonzero characteristic roots are simple. Each of these operators, taken as $Y_1$, can be replaced in the described way by several operators of the form $(X_{\alpha}, X_{\alpha'})$ for each of which the nonzero characteristic roots are simple and are rationally related. Thus, as the operators generating the group $\mathcal{G}$, one can take the independent operators of this kind $Z_1, Z_2, \dots, Z_u$. Now, in order to find the finite equations of the subgroup $\mathcal{H}$, it is sufficient to find the invariants of degree $m$ of the group $\mathcal{G}$. Such will be the solutions of the system
\begin{equation}\label{eq19} 
Z_i + Z_i' + \cdots + Z_i^{(m-1)} = 0 \quad (i=1,2,\dots,u), \tag{1.9}
\end{equation}
where the space $(x_1, x_2, \dots, x_r)$ is the representation space of the adjoint group, and the spaces $(x_1^{(i)}, x_2^{(i)}, \dots, x_r^{(i)})$, where $(i=0,1,\dots,m-1)$, together contain the parametric group. The solutions of each individual equation of the system \ref{eq19} are integrals of one of the systems of linear homogeneous differential equations
\begin{equation}\label{eq110}
\frac{dx_1}{\sum_\nu c_{i\nu}^1 x_\nu} = \frac{dx_2}{\sum_\nu c_{i\nu}^2 x_\nu} = \cdots = \frac{dx_r^{(m-1)}}{\sum_\nu c_{i\nu}^r x_\nu^{(m-1)}} \quad\quad (i = 1,2,\dots,u). \tag{1.10} 
\end{equation}
Each of these systems has simple elementary divisors, and therefore its integrals can be written in the form
\begin{equation*}
z_1 = c_1 e^{\omega_1 t} , z_2 = c_2 e^{\omega_2 t}, \dots, z_r^{(m-1)} = c_r^{(m-1)} e^{\omega_r t}
\end{equation*}  
where the $z_i$ represent some linear functions of $x_i$. Since the ratios of the quantities $\omega_i$ are rational, algebraic integrals are obtained from this by eliminating $t$. To obtain solutions of the whole system, one needs to construct such functions that can be represented as functions of the integrals of each of the systems \ref{eq110} for $i=1,2,\dots,u$, which is a purely algebraic problem. Thus we have proven

\begin{theorem}
Every semisimple $s$-group has an algebraic representation in $s$-dimensional space (\cite[p. 133]{cartan1894structure}).
\end{theorem}

Killing and Cartan (\cite[p. 151]{cartan1894structure}) determined all the possible types of simple groups and calculated for each of them the index $s$ of the largest subgroup. Their results are presented as follows\footnote{Translator note: In the original text, Chebotarev denotes the exceptional cases $E_6$, $E_7$, $E_8$, $F_4$, and $G_2$  by $G_{78}$, $G_{133}$, $G_{248}$, $G_{52}$, and $G_{14}$, respectively. We chose to use the now standard notation for these simple Lie algebras in the main text.}:

\begin{enumerate}
\item Unimodular linear homogeneous groups in $n$ variables (type $A$); $r=n^2 - 1$, $s=n-1$.
\item Orthogonal groups in $n$ variables (type $B$ for odd $n$ and type $D$ for even $n$); $r=\frac{n(n-1)}{2}$, $s=n-1$. 
\item Complex groups in $n$ variables ($n$ even), i.e. linear homogeneous groups that leave the Pfaffian
\begin{equation*}
(x_1dx_2 - x_2dx_1) + (x_3dx_4 - x_4dx_3) + \cdots + (x_{n-1}dx_n - x_ndx_{n-1}).
\end{equation*}
invariant under transformation (type $C$); $r=\frac{n(n+1)}{2}, s=n$.
\item Type $E$: \\
1) The group $E_6$; $r=78$, $s=16$; \\
2) The group $E_7$; $r=133$, $s=27$; \\  
3) The group $E_8$; $r=248$, $s=57$.
\item Type $F$: $F_4$; $r=52$, $s=15$.
\item Type $G$: $G_2$; $r=14$, $s=5$.
\end{enumerate}

\section{Finite and continuous groups} 

The connection between finite and continuous groups has been studied for a long time, when the following problem was posed, which led to very remarkable results (\cite{jordan1878memoire} and \cite{speiser1927theorie}): to determine all finite subgroups of a given continuous group. 

We will deal with the following reversal of this problem: to determine all ctg's isomorphic to a given finite group.\footnote{Translator's note: Chebotarev probably intended to write ``determine all ctg's with a subgroup isomorphic to a given finite group" --- in the main text we give the direct translation.}

Since with this formulation of the problem we are forced to consider an infinite set of solutions uninteresting to us (for example, any group containing any solution will also be a solution), we will refine our problem by calling it the problem of finding \textit{enveloping groups}\footnote{In previous works I used the name ``Einkleidungsgruppe" (here the author uses the term ``Einbettungsruppe" - Ed.)} (for short: eg) and define this concept as follows:

A ctg $\Gamma$ is called an eg of a finite group $\mathfrak{G}$ if: 
\begin{enumerate}
    \item $\Gamma$ contains a divisor isomorphic to $\mathfrak{G}$, 
    \item none of the proper continuous subgroups of the group $\Gamma$ has property 1,
    \item none of the proper continuous factor groups of the group $\Gamma$ has property 1.
\end{enumerate}

In most cases when we have to look for an eg for a finite group $\mathfrak{G}$, the latter will be a simple group. Then the following is true: 

\begin{theorem} 
Every eg of a simple group $\mathfrak{G}$ is also simple.
\end{theorem}

\begin{proof}
Suppose the eg $\Gamma$ of the group $\mathfrak{G}$ contains a proper continuous normal divisor $\Gamma_1$. Say the subgroup $G$ of $\Gamma$, isomorphic to $\mathfrak{G}$, has intersection $G_1$ with $\Gamma_1$. Since each subgroup $S^{-1}G_1S$ conjugate to the group $G_1$ is contained in the group $S^{-1}\Gamma_1S = \Gamma_1$, $G_1$ is a normal divisor of the group $G$. 

The latter is simple, and therefore either $G = G_1$ or $G_1 = 1$. In the first case, we arrive at a contradiction with condition 2. In the second case, among the elements $U\Gamma_1$ of the quotient group $\Gamma/\Gamma_1$ (where $U$ runs through the elements of the group $\Gamma$) are contained all the elements $SG_1$, where $S$ runs through the whole set of elements of the group $G$. The elements $SG_1$ are all distinct from each other and therefore form a group isomorphic to $G$, which contradicts condition 3.
\end{proof}

When we look for an eg of a given finite group $\mathfrak{G}$, it is natural to define for each eg one of its representations, for example, as a linear homogeneous group. But here a certain difficulty arises, since all representations of the same ctg are only locally isomorphic to each other. Therefore, it may well happen that one of two representations contains the group $\mathfrak{G}$, while the other does not. Consequently, we must delve a little deeper into the theory of locally isomorphic ctg's developed by O. Schreier \cite{schreier1927verwandtschaft}. 

If one considers the elements of each ctg as points of the corresponding ``group space", then the local isomorphism of two ctg means the mutual correspondence of neighborhoods of the identity elements of these groups. This one-to-one correspondence can be analytically continued indefinitely, since the neighborhoods of any corresponding points lying inside these corresponding neighborhoods can also be mutually mapped onto each other. However, from this one can conclude that both groups are isomorphic (as a whole) only if each of the corresponding group spaces are simply connected, i.e. if in each of them every closed path can be continuously contracted to a point. If the group space is not simply connected, then, as Schreier proved (cited above), one can associate with it a certain simply connected covering space. The group (called the \textit{covering group}) corresponding to this space has the property that every locally isomorphic group is isomorphic to its quotient group relative to some normal divisor $\Delta$ whose elements are discrete (this means that for each element of $\Delta$ one can find such a neighborhood inside the group space that does not contain other elements of $\Delta$). 

\begin{theorem}[Schreier] 
If two ctg's $\Gamma$ and $\Gamma_1$ are locally isomorphic and if the group $\Gamma_1$ is isomorphic to the group $\Gamma/\Delta$, then $\Delta$ lies in the center of the group $\Gamma$.
\end{theorem}

\begin{proof}
Let $U$ be an arbitrary transformation of the group $\Gamma$. Continuously vary the parameter values of the group $\Gamma$ so that the corresponding transformation continuously changes from the identity transformation $I$ to $U$. Since $\Delta$ is a normal divisor of the group $\Gamma$, all transformations $U^{-1}DU$ (where $D$ belongs to $\Delta$) are contained in $\Delta$. On the other hand, $I^{-1}DI=D$. With continuous variation of the parameters of the transformation $U$, the transformation $U^{-1}DU$ also changes continuously and therefore remains in some neighborhood of the transformation $D$ not containing any other transformations from $\Delta$. Consequently, $U^{-1}DU=D$, and the theorem is proved. 
\end{proof}

Let $\Gamma_1$ be an eg of the finite group $\mathfrak{G}$ and $\Gamma$ the covering group of the group $\Gamma_1$. $\Gamma$ may not contain any subgroup isomorphic to $\mathfrak{G}$, since, when $\Gamma$ is mapped onto $\Gamma_1$, such a group $G$ can pass into a subgroup $G_1$ isomorphic with the group $G$ with only some quotient group $G/H$ is isomorphic to $\mathfrak{G}$. Here $H$ is a divisor of $\Delta$ and therefore is contained in the center of the group $G$. 

In particular, if the group $\mathfrak{G}$ is simple, then the commutator group $G'$ of the group $G$ also contains the group isomorphic to $\mathfrak{G}$ as a factor group and is contained in $\Gamma$ as a divisor. Let us prove that $G'$ coincides with its commutator group $G''$. Since the factor group $G/\Delta$ is isomorphic to the group $\mathfrak{G}$, i.e. simple, and therefore coincides with its commutator group, any element $s$ of the group $G/\Delta$ can be represented as 
$$\prod_{i,j} s_is_js_i^{-1}s_j^{-1}.$$
It follows from this that any element $S$ of the group $G$ can be represented in the form
\begin{equation}\label{eq21}
    S = D\prod_{i,j}S_iS_jS_i^{-1}S_j^{-1}
    \tag{2.1}
\end{equation}
where $D$ is an element of the group $\Delta$. In order for all elements of $G'$ to be contained in $G''$, it is sufficient to show this for elements of the form $S_1S_2S_1^{-1}S_2^{-1}$. If 
\begin{align*}
S_1 &= D_1\prod_{i,j}S_{1i}S_{1j}S_{1i}^{-1}S_{1j}^{-1}, \\  
S_2 &= D_2\prod_{i,j}S_{2i}S_{2j}S_{2i}^{-1}S_{2j}^{-1},
\end{align*}
then 
\begin{align*}
S_1S_2S_1^{-1}S_2^{-1} &= \left(D_1\prod_{i,j}S_{1i}S_{1j}S_{1i}^{-1}S_{1j}^{-1}\right)\left(D_2\prod_{i,j}S_{2i}S_{2j}S_{2i}^{-1}S_{2j}^{-1}\right) \\
&\quad \cdot \left(D_1\prod_{i,j}S_{1i}S_{1j}S_{1i}^{-1}S_{1j}^{-1}\right)^{-1}\left(D_2\prod_{i,j}S_{2i}S_{2j}S_{2i}^{-1}S_{2j}^{-1}\right)^{-1}
\end{align*}
and since the elements $D_1, D_2$ commute with all $S$,
\begin{align*} 
S_1S_2S_1^{-1}S_2^{-1} &= \prod_{i,j}\left(S_{1i}S_{1j}S_{1i}^{-1}S_{1j}^{-1}\right)\prod_{i,j}\left(S_{2i}S_{2j}S_{2i}^{-1}S_{2j}^{-1}\right) \\
&\quad \cdot \prod_{i,j}\left(S_{1i}S_{1j}S_{1i}^{-1}S_{1j}^{-1}\right)^{-1}\prod_{i,j}\left(S_{2i}S_{2j}S_{2i}^{-1}S_{2j}^{-1}\right)^{-1}.
\end{align*}
It follows immediately from this that $S_1S_2S_1^{-1}S_2^{-1}$ is contained in $G''$, which was to be proved. 

Thus it is proved that the covering group of the group $\Gamma_1$, and therefore any ctg (locally) isomorphic to $\Gamma_1$, has a divisor $G'$ whose quotient group relative to the center is isomorphic to the group $\mathfrak{G}$ and whose commutator group coincides with the group $G'$ itself.

Such groups can be considered a special case of the representation groups $\mathfrak{K}$ of the group $\mathfrak{G}$ examined by I. Schur \cite{schur1904darstellung}\footnote{Translator's note: Chebotarev is recalling what we now call the theory of projective representations and Schur multipliers.}, which he characterized by the following properties:
\begin{enumerate}
    \item $\mathfrak{K}$ contains a subgroup $\mathfrak{M}$ consisting of invariant elements of the group $\mathfrak{K}$, such that the quotient group $\mathfrak{K}/\mathfrak{M}$ is isomorphic to the group $G$;
    \item The commutator group of $\mathfrak{K}$ contains all elements of the group $\mathfrak{M}$;
    \item There is no group that has properties 1 and 2 and whose order exceeds the order of the group $\mathfrak{K}$.
\end{enumerate}
I. Schur shows that every group $\mathfrak{G}$ has only a finite number of non-isomorphic representation groups and gives a method for constructing all representation groups. In doing so, he proves that if the group $\mathfrak{G}$ coincides with its commutator group, which certainly takes place if the group $\mathfrak{G}$ is simple, then there exists only one representation group of the group $\mathfrak{G}$ (\cite[p. 38, IV]{schur1904darstellung}).

\begin{example}
The alternating group of degree six can be represented as a group of fractional linear substitutions in two variables. However, it cannot be represented as a homogeneous linear group in three variables. The latter takes place for the first time for its representation group whose order is equal to 3,360 \cite{wiman1896ueber}.
\end{example}

In order to construct an eg of a certain finite group, note that any representation $\Gamma$ of any simple ctg can be represented as a linear homogeneous group. In doing so, the subgroup $G$ of this group isomorphic to the group $\mathfrak{G}$ will be represented as a finite linear homogeneous group. Let us assume that the group $\mathfrak{G}$ is reducible, namely that its matrices by means of the linear transformation $T$ are reduced to the form 
$$S=\begin{pmatrix} A & 0 \\ 0 & B \end{pmatrix}.$$
It does not yet follow from this that the group $T^{-1}\Gamma T$ is reducible. However, if one takes any infinitesimal operator $X$ of the group $T^{-1}\Gamma T$ containing the substitution $S$, then $X$ (considered as a matrix) will commute with $S$. If all the characteristic roots of the matrix $S$ were simple, then the same linear transformation would reduce the matrices $X$ and $S$ to diagonal form simultaneously. Therefore, $X$ has the form 
$$X = \begin{pmatrix} M & 0 \\ 0 & N \end{pmatrix}$$ 
On the other hand, if $S$ runs through the system of generating substitutions of the group $\mathfrak{G}$ (i.e. those whose composition gives the whole group $\mathfrak{G}$), then the system of infinitesimal operators $X_i$ corresponding to them cannot be contained in any proper continuous subgroup of the group $T^{-1}\Gamma T$. Indeed, otherwise, by condition 2 of the definition of an eg, the group $\Gamma$ would not be an eg of the group $\mathfrak{G}$. In other words, the operators $X_i$ and the operators $(X_i, X_j)$ generated by them constitute a complete system of infinitesimal operators of the group $T^{-1}\Gamma T$. It follows from this that the group $\Gamma$ can be irreducible only when each system of generating substitutions of the group $\mathfrak{G}$ contains substitutions with multiple characteristic roots. This circumstance must necessarily take place if some irreducible parts of the group $\mathfrak{G}$ are similar to each other. It seems very likely that characteristic roots can be multiple only in this case. 

The general problem of enveloping the group $\mathfrak{G}$ can be formulated as follows. Let $I, A, B, \dots$ be the set of linear substitutions corresponding to the elements of the group $\mathfrak{G}$. Let us introduce the notation 
$$X^A = A^{-1}X A$$
and let $X_1, X_2, \dots, X_r$ be the sought infinitesimal operators of the eg $\Gamma$ (which we consider to be a linear homogeneous group) of the group $\mathfrak{G}$. Then for each $A$ the equations
$$X_i^A = \sum_{\nu =1}^r a^A_{i\nu}X_\nu \quad (i = 1,2,\dots,r),$$
must hold, where $a^A_{i\nu}$ are unknown constants. The correspondence $A \mapsto (a_{i\nu}^A)$ gives a new linear homogeneous representation of the group $\mathfrak{G}$. We will call it a representation ``reducible in the narrow sense" if there exist such linear combinations $Y_1, Y_2, \dots, Y_m$ of the operators $X_i$ that: 
\begin{enumerate}
    \item  form a subgroup of the group $\Gamma$ (i.e. $(Y_i, Y_j)$ are linearly expressed in terms of $Y_i$), 
    \item undergo only linear substitutions by means of the transformations $A$.
\end{enumerate}
It is obvious that the ctg $\Gamma$ corresponding to the representation reducible in the narrow sense is not an eg of the group $\mathfrak{G}$, since it has a proper continuous subgroup containing the group isomorphic to $\mathfrak{G}$ as a divisor. Thus, our problem is reduced to finding systems $X_1, X_2, \dots, X_r$ that are irreducible in the narrow sense. 

It is very difficult to solve this general problem. I cannot even say so far whether each finite group corresponds to a finite or infinite number of eg. But if we set ourselves the task of determining for a given finite simple group $\mathfrak{G}$ only those eg that have a representation in a space of the smallest possible dimension $s$, then this problem can be solved relatively easily. To do this, one must find all irreducible representations of both the group $\mathfrak{G}$ and its representation groups. Then, choosing a proper representation $\mathfrak{g}$ of the smallest possible number of dimensions $f$, one must determine whether $\mathfrak{g}$ is complex or real. In the first case, the full unimodular linear homogeneous group of degree $f$ will be the sought eg if the locally isomorphic group of fractional linear substitutions in $(f-1)$ variables actually contains a divisor isomorphic to $\mathfrak{G}$. The latter takes place when $\mathfrak{g}$ is a representation of the group $\mathfrak{G}$, and not of its representation group. If the group $\mathfrak{g}$ is real, then it is reduced to a real orthogonal group (\cite[p. 107, Theorem 100]{speiser1927theorie}) and therefore will be a divisor of the orthogonal group of degree $f$ which is an $(f-2)$-group. This method is, up to a finite number of exceptions, the most general method for constructing eg having representations in a space of the smallest possible number of dimensions. Indeed, apart from a finite number of exceptions, every simple ctg can be represented either as a full unimodular linear homogeneous group, as a complex group, or as a full orthogonal group. But the second case is of no interest to us, since the transition from a full unimodular linear homogeneous group to a complex group does not give a decrease in the number of dimensions of the smallest representation space. The subgroup isomorphic to $\mathfrak{G}$ enters these representations either as a linear homogeneous or as an orthogonal group. It will be irreducible, except in the case when the smallest value of $f$ corresponds to a representation group, so that its enveloping group does not contain subgroups isomorphic to $\mathfrak{G}$. The value of $f$ corresponding to the group $\mathfrak{G}$ itself will in this case be greater than the sum of two (or more) values of $f$ corresponding to the representation groups. Thus, depending on whether the representation corresponding to the smallest value of $f$ is complex or real, we will have either $s=f-1$ or $s=f-2$.

In particular, if $\mathfrak{G}$ is an alternating group of $n$ variables and $n \geq 8$, then, as A. Wiman proved \cite{wiman1899ueber}, neither the group $\mathfrak{G}$ itself nor its representation group can be represented as linear homogeneous groups of degree less than $(n-1)$. On the other hand, the representation of the group $\mathfrak{G}$ as an alternating group of substitutions in $n$ variables decomposes into the identical representation and a real irreducible representation of degree $(n-1)$. Thus, here $s = n-3$.

One question remains to be solved: is there a ctg $\Gamma_1$ that can be represented in a space of fewer dimensions than a given simple ctg $\Gamma$ isomorphic to one of the factor groups of $\Gamma_1$? Such a case is impossible. Indeed, E. E. Levi \cite{levi1905sulla} showed that every ctg $\Gamma_1$ having a simple factor group $\Gamma$ contains a subgroup isomorphic to $\Gamma$. Therefore, if the group $\Gamma_1$ is representable in an $s$-dimensional space, the same is obviously true for its subgroups and hence for the factor group $\Gamma$.

\section{Solution of Klein's resolvent problem}

The first problem formulated in the introduction is resolved by the following theorem:

\begin{theorem}\label{thm8}
If the Galois group $\mathfrak{G}$ of the equation 
\begin{equation}\label{eq31}
f(x) = x^n + a_1x^{n-1} + \cdots + a_{n-1}x + a_n = 0 \tag{3.1}
\end{equation}
is simple and, considered as an abstract group, is a divisor of some ctg $\Gamma$ and if $\Gamma$ is an $s$-group, then the equation \ref{eq31} has an $s$-parametric resolvent. 
\end{theorem}

\begin{remark}
There is no doubt that this theorem is valid for any finite group $\mathfrak{G}$. I assume the group $\mathfrak{G}$ to be simple only because at present the possibility of rationally performing all the necessary operations seems directly obvious only for simple groups.
\end{remark}

First we will prove the following theorem:

\begin{theorem}
Equation \ref{eq31} has an $s$-parametric resolvent if and only if its Galois resolvent
\begin{equation}\label{eq32}
\Phi_n(y) = 0 \tag{3.2}
\end{equation}
with arbitrarily variable coefficients has one.
\end{theorem}

\begin{proof}
Let equation \ref{eq31} have an $s$-parametric resolvent
\begin{equation}\label{eq33}
f_1(y) = 0, \tag{3.3}  
\end{equation}
whose coefficients are rational functions of $a_1, a_2, \dots, a_n; \Phi$, where $\Phi$ is some function of the roots $x_1, x_2, \dots, x_n$ of equation \ref{eq31} belonging to the group $\mathfrak{G}$. Let, in addition, the given normal equation \ref{eq32} with arbitrary variable coefficients be given, whose group $\mathfrak{G}$ is simply transitive. Denote by $\overline{\Phi}$ the function of the roots of equation \ref{eq32} belonging to the group $\mathfrak{G}$. To construct an $s$-parametric resolvent of equation \ref{eq32}, choose a function $x(\xi)$, belonging to the same group $\mathfrak{H}$ as the root $x_1$ of equation \ref{eq31} within the field $K(x_1, x_2, \dots, x_n)$. We will consider the coefficients of the $n$th degree equation satisfied by $x(\xi)$ as partial values of the coefficients $a_i$ of equation \ref{eq31}. In the same way, express $\Phi$ through $A_i, \overline{\Phi}$. Then the coefficients of both the resolvent \ref{eq33} and the transition formulas from \ref{eq31} to \ref{eq33} will be expressed through $A_i, \overline{\Phi}$. Now construct the Galois resolvent for equation \ref{eq32}
\begin{equation} \label{eq34}
F_1(\eta) = 0. \tag{3.4}
\end{equation}
It can be considered an $s$-parametric resolvent of equation \ref{eq32}, since any root $\xi$ of equation \ref{eq32} is rationally expressed through the roots of equation \ref{eq31}, which in turn are rationally expressed through the roots of equation \ref{eq33} and the quantities $a_1, a_2, \dots, a_n; \Phi$, i.e. through one of the roots of equation \ref{eq34} and the quantities $A_1, A_2, \dots, A_m; \overline{\Phi}$.

Conversely, if equation \ref{eq32} has an $s$-parametric resolvent \ref{eq34} and the root $x_1$ of equation \ref{eq31} belongs to the group $\mathfrak{H}$, then the coefficients $A_1, A_2, \dots, A_m; \overline{\Phi}$ of equation \ref{eq32} are rationally expressed through the coefficients $a_1, a_2, \dots, a_n; \Phi$ of equation \ref{eq31}. In this case, the roots of equation \ref{eq33} belong to the group $\mathfrak{H}$ within the field $K(\eta)$, where $\eta$ is one of the roots of equation \ref{eq34}, and therefore equation \ref{eq33} can be considered an $s$-parametric resolvent of equation \ref{eq31}.
\end{proof}

\begin{proof}[\proofname\ of Theorem \ref{thm8}]
Let the roots of equation \ref{eq31}, which we consider as independent variables, be given by
$$x_1, x_2, \dots, x_n.$$ 
In addition, let a ctg $\overline{\Gamma}$ be given that contains as a divisor the group $\overline{\mathfrak{G}}$ isomorphic to the Galois group $\mathfrak{G}$ of equation \ref{eq31} (i.e. the group to which the function $\Phi(x_1, x_2, \dots, x_n)$ associated to the field of coefficients belongs). Let the group $\overline{\Gamma}$ be an $s$-group, namely having a representation in the space $(z_1, z_2, \dots, z_s)$.

We can consider the group $\overline{\Gamma}$ to be simple. Indeed, if the group $\overline{\Gamma}$ has a proper normal divisor $\overline{\Gamma}_1$, then either $\overline{\Gamma}_1$ or $\Gamma/\Gamma_1$ has a divisor isomorphic to the group $\mathfrak{G}$. In this case, by E. E. Levi's theorem \cite{levi1905sulla}, these groups are $s$-groups. Continuing this process, we will eventually arrive at a simple eg of the group $\mathfrak{G}$, which will also be an $s$-group.

Denote by $\Gamma$ the representation of our $s$-group either as a full unimodular linear homogeneous group or as a full orthogonal group, whose $m$th power contains the representation $\overline{\Gamma}$ (in the first case $m=1$). In both cases, it is easy to see that the coordinates $z_1, z_2, \dots, z_s$ of the representation space of the group $\overline{\Gamma}$ are rationally expressed in terms of the coordinates of the corresponding space of the group $\Gamma^m$. 

We will consider the group of substitutions $\mathfrak{G}$ as a linear homogeneous group in the variables $x_1, x_2, \dots, x_n$ and subject these variables to such a linear transformation that the group $\mathfrak{G}$ is transformed into a fully reduced form. Then consider those components of this representation that are divisors of the representation $\Gamma^m$. In other words, we seek such linear functions $y_1, y_2, \dots, y_u$ of the variables $x_1, x_2, \dots, x_n$ that will undergo substitutions corresponding to the above mentioned irreducible components of the group $\mathfrak{G}$ if the variables $x_1, x_2, \dots, x_n$ are subjected to substitutions of the group $\mathfrak{G}$. If equation \ref{eq31} is normal, then the representation of $\mathfrak{G}$ is regular and therefore contains as constituent parts all irreducible linear representations of this group, each as many times as its degree (\cite[p. 119, note 14]{speiser1927theorie}). If these irreducible representations enter as divisors in $\Gamma^m$ with higher multiplicities than in $\mathfrak{G}$, then we will add to the system $x_1, x_2, \dots, x_n$ additional parallel systems of $n$ variables and denote the number of these systems again by $m$
\begin{equation}\label{eq35}
    \begin{matrix}
        y_1, y_2, \dots, y_s, \\
y_1^{(1)}, y_2^{(1)}, \dots, y_s^{(1)}, \\
 \vdots \\
y_1^{(m-1)}, y_2^{(m-1)}, \dots, y_s^{(m-1)}
    \end{matrix}\tag{3.5}
\end{equation}
We can consider the variables \ref{eq35} as coordinates of the representation space of the group $\Gamma^m$. Subjecting the variables $x_1, x_2, \dots, x_n$ to substitutions of the Galois group of equation \ref{eq31} and subjecting all systems of variables \ref{eq35} to the same substitutions, we obtain a divisor contained in this representation that is isomorphic to the group $\mathfrak{G}$. Since $\overline{\Gamma} \preceq\Gamma^m$; we can find such rational functions 
\begin{equation}\label{eq36}
    z_1,z_2,\dots, z_s 
    \tag{3.6}
\end{equation}
of the variables \ref{eq35} where transformations of the group $\Gamma^m$, and hence by substitutions $\mathfrak{G}^m$ as well, transform into functions of $z_1, z_2, \dots, z_s$. 

The first system \ref{eq35} consists of certain homogeneous linear functions of the roots $x_1, x_2, \dots, x_n$ of equation \ref{eq31}. We will consider the other systems \ref{eq35} as the same functions of new variables $x_1^{(i)}, x_2^{(i)}, \dots, x_n^{(i)}$, so that we get the following new system of variables
\begin{equation}\label{eq37}
    \begin{matrix}
x_1, x_2, \dots, x_n, \\
x_1^{(1)}, x_2^{(1)}, \dots, x_n^{(1)},\\
\vdots \\
x_1^{(m-1)}, x_2^{(m-1)}, \dots, x_n^{(m-1)}, 
\end{matrix}
\tag{3.7}
\end{equation}

Let us transform the system of variables \ref{eq37} as follows:
\begin{equation}\label{eq38}
    \begin{matrix}
x_1^{(i)} = \alpha_0^{(i)} + \alpha_1^{(i)}x_1 + \cdots + \alpha_{n-1}^{(i)}x_{n-1}, \\
x_2^{(i)} = \beta_0^{(i)} + \beta_1^{(i)}x_1 + \cdots + \beta_{n-1}^{(i)}x_{n-1}, \\
\vdots \\
x_n^{(i)} = \gamma_0^{(i)} + \gamma_1^{(i)}x_1 + \cdots + \gamma_{n-1}^{(i)}x_{n-1}
\end{matrix}
\tag{3.8}
\end{equation}
where $(i = 1,2,\dots, m-1)$, and take as new variables the quantities 
\begin{equation}\label{eq39}
    x_1,x_2,\dots,x_n;\alpha_0^{(i)}, \alpha_1^{(i)}, \dots, \alpha_{n-1}^{(i)}, \quad (i=1,2,\dots,m-1). \tag{3.9}
\end{equation}
Obviously, the transition from \ref{eq37} to \ref{eq39} is carried out by a reversible transformation. If we now perform substitutions of the group $\mathfrak{G}^m$ on the system \ref{eq39}, the variables 
$$\alpha_0^{(i)}, \alpha_1^{(i)}, \dots, \alpha_{n-1}^{(i)} \quad (i=1,2,\dots,m-1)$$
do not change.

Substituting into expression \ref{eq36} instead of $y_j^{(i)}$ their values calculated using formulas $\ref{eq38}$, we obtain a system of $s$ functions 
\begin{equation}\label{eq310}
\overline{Z}_1, \overline{Z}_2, \dots, \overline{Z}_s \tag{3.10}
\end{equation}
of the variables \ref{eq39}, and these functions will be transformed into their own functions if substitutions are made on the variables $x_1, x_2, \dots, x_n$ from the group $\mathfrak{G}$, leaving the quantities $\alpha_j^{(i)}$ invariant. Therefore, if in one of the functions \ref{eq310}, say $\overline{Z}_1$, we assign any rational values to the quantities $\alpha_j^{(i)}$ and then subject the variables $x_1, x_2, \dots, x_n$ to substitutions from the group $\mathfrak{G}$, we will obtain a system of functions
\begin{equation}\label{eq311}
\overline{Z}, \overline{Z}^{S_2}, \dots, \overline{Z}^{S_n},\tag{3.11}
\end{equation}
functionally dependent on $\overline{Z}_1, \overline{Z}_2, \dots, \overline{Z}_s$, i.e. containing only $s$ variable parameters. The elementary symmetric functions of the quantities \ref{eq311} are rationally expressed through $a_1, a_2, \dots, a_n; \Phi$ and also contain only $s$ variable parameters (in other words, contain only $s$ functionally independent quantities). To make sure that we arrive at an $s$-parametric resolvent of equation $\ref{eq31}$, it is sufficient to show that with a suitable choice of functions $z_i$ and values of $\alpha_j^{(i)}$, the quantities \ref{eq311} are distinct from each other. 

Suppose this is not the case. Then in the group $\mathfrak{G}$ there exists a substitution $S$ where all the differences $\overline{Z}_i^S - \overline{Z}_i$ vanish for all values of $\alpha_j^{(i)}$. Moreover: each of the functions $Z_{i}^U$, where $U$ is an arbitrary transformation of the ctg $\overline{\Gamma}^m$, satisfies the equation 
$$Z_{i}^{US} - Z_{i}^U = 0.$$
Applying to this equation the transformation $U^{-1}$, we obtain
$$Z_i^{USU^{-1}} - Z_i = 0.$$
But $S$ cannot commute with all transformations of the group $\overline{\Gamma}^m$. Consequently, the set of transformations $USU^{-1}$ represents such a continuous family of transformations that contains $S$ and other transformations distinct from $S$. The composite of this family (i.e. the set of transformations obtained by composing transformations of the form $USU^{-1}$) therefore represents a continuous group that leaves invariant all functions $Z_1, Z_2, \dots, Z_s$. This contradicts the assumption that the representation of the group $\overline{\Gamma}^m$ in the space \ref{eq36} is proper. If one of the differences $Z_i^S - Z_i$ does not vanish identically, then it is possible to choose rational values of the quantities $\alpha_j^{(i)}$ so that all the quantities \ref{eq311} are distinct from each other and therefore represent all the roots of the $s$-parametric resolvent of equation \ref{eq31}. 

These arguments must be modified somewhat in the case when not the group $\mathfrak{G}$ itself but one of its representation groups, say $\mathfrak{H}$, is a divisor of the linear homogeneous ctg $\overline{\Gamma}$ whose degree $m$ group $\overline{\Gamma}^m$ contains the $s$-dimensional representation. Then we solve the resolvent problem for the group $\mathfrak{G}$ by constructing the equation 
\begin{equation}\label{eq312}
    F(\xi) = 0, \tag{3.12}
\end{equation}
whose group is $\mathfrak{H}$ and whose roots are rationally expressed through the roots of equation \ref{eq31}. To do this, it is sufficient to consider the equation 
$$f(\zeta^d) = 0,$$
where $d$ is the order of the group $\mathfrak{F}$ (here $\mathfrak{G} \cong \mathfrak{H}/\mathfrak{F}$), and construct such a subfield of the field $K(\zeta)$ whose group is isomorphic to $\mathfrak{G}$. The further reasoning can be repeated without change, since the eg of the group $\mathfrak{G}$ is simple, although $\mathfrak{G}$ itself is not. The roots of the resolvent constructed in this way belong in the field $K(\zeta)$ not to the identity group, but to the group $\mathfrak{F}$. But the roots of equation \ref{eq31} belong to the same group as well. Consequently, they are rationally expressed through the roots of the resolvent and the quantities $a_1, a_2, \dots, a_n; \Phi$.
\end{proof}

\section{Inversion of the main theorem}

Theorem \ref{thm8}, which due to its significance for the whole theory, is the main theorem. It admits the following inversion:

\begin{theorem}\label{thm10}
If equation \ref{eq31} has an $s$-parametric resolvent, then its group $\mathfrak{G}$ has as eg some $s$-group. 
\end{theorem}

\begin{proof}
Let the roots of the $s$-parametric resolvent of equation \ref{eq31}, which we consider normal, be
\begin{equation}\label{eq41}
Z_1, Z_2, \dots, Z_s. \tag{4.1}
\end{equation}
The quantities \ref{eq41} are such functions of the roots $x_1, x_2, \dots, x_n$ of equation \ref{eq31} that each of them belongs to the identity group. If we perform substitutions from the group $\mathfrak{G}$ on $x_1, x_2, \dots, x_n$, then the quantities \ref{eq41} will undergo the known substitutions that form the group $\overline{\mathfrak{G}}$, isomorphic to $\mathfrak{G}$. The groups $\overline{\mathfrak{G}}$ and $\mathfrak{G}$ differ from each other only in the notation of the variables. 

Let us now consider the quantities \ref{eq41} as independent variables and form the group $\overline{\mathfrak{G}}$ by means of an arbitrary ctg $\Gamma$. We can use a linear homogeneous group in the variables \ref{eq41} for this. It is even better if, by means of a linear transformation of the quantities $Z_i$, we represent the group $\overline{\mathfrak{G}}$ in a fully reduced form and then form one of its irreducible parts by means of a linear ctg.

Then we utilize the fact that the quantities \ref{eq41} represent the roots of the $s$-parametric resolvent, i.e. that there exist exactly $s$ functionally independent quantities between them. We will consider $x_1, x_2, \dots, x_n$ to be the coordinates of a special space\footnote{Translator's note: The phrase ``special space" just means that $\mathfrak{R}$ is given the quotient topology.} $\mathfrak{R}$ in which we will consider two points $(x_1, x_2, \dots, x_n)$ and $(x_1', x_2', \dots, x_n')$ to coincide if 
\begin{equation}\label{eq42}
    Z_i(x_1, x_2, \dots, x_n) = Z(x_1', x_2', \dots, x_n') \quad (i=1,2,\dots,s). \tag{4.2}
\end{equation}
Obviously, the space $\mathfrak{R}$ is $s$-dimensional. 

Let now the transformations of the ctg $\Gamma$ be
\begin{equation}\label{eq43}
Z_i' = f_i(Z_1, Z_2, \dots, Z_s; a_1,a_2,\dots,a_r) \quad (i=1,2,\dots,s). \tag{4.3}
\end{equation}
If we substitute here instead of $Z_i$ and $Z_i'$ their values computed using the formulas $Z_i(x_1, x_2, \dots, x_n)$ and $Z_i(x_1', x_2', \dots, x_n')$, then the obtained equations will define a transformation ``induced in the space $\mathfrak{R}$" by the ctg $\Gamma$, which will contain the group $\mathfrak{G}$ as a divisor. Indeed, the collection of functions $Z_i(x_1, x_2, \dots, x_n)$ belong to the identity group and therefore those substitutions over the variables $x_1, x_2, \dots, x_n$ which are contained in the group $\mathfrak{G}$ actually produce substitutions over the functions $Z_i$ and hence, according to the definition of $\mathfrak{R}$, correspond to distinct transformations of the space $\mathfrak{R}$. The theorem is proved.
\end{proof}

\begin{remark} 
For the proof of Theorem \ref{thm10}, it is essential that each system of parameter values of the ctg $\Gamma$ uniquely determines the point to which a given point of the space $\mathfrak{R}$ goes. By virtue of this, the group $\Gamma$ falls under Schreier's (\cite[note 16]{schreier1927verwandtschaft}) topological definition, so that his whole theory applies to it. If this condition is not met, one can, for example, encounter a phenomenon where the non-cyclic monodromy group of some algebraic function has a one-parameter ctg as its eg. For example: the one-parameter ctg of the space $(x, y, z)$ defined by the equations 
\begin{align*}
    x + y + z = c_1 t,\\
    x^2 + y^2 + z^2 = c_2t,\\
    x^3 + y^3 + z^3 = c_3t,
\end{align*}
contains as a divisor the symmetric group of substitutions of the third degree, which is not even abelian.
\end{remark}

Together with the above quoted studies of A. Wiman (\cite[note 21]{wiman1899ueber}), Theorem \ref{thm10} makes it very plausible to assume that for $n \geq 8$
$$n - s \leq 3.$$
For Hilbert's problem, as A. Wiman recently showed \cite{wiman1927anwendung},
$n - s \geq 5$.

%
%

\newpage
\bibliographystyle{alpha} 
\bibliography{chebotarevKH}

\end{document}